\documentclass[12pt,preprint]{article}

\usepackage{amsmath, latexsym, amsfonts, amssymb, amsthm}
\usepackage{graphicx, color,hyperref,dsfont,epsfig,caption,wrapfig,subfig}
\usepackage{pstricks}
\usepackage{movie15}
\hypersetup{
    linktoc=page,
    linkcolor=red,          
    citecolor=blue,        
    filecolor=blue,      
    urlcolor=cyan,
    colorlinks=true           
}

\numberwithin{equation}{section}

\newtheorem{theorem}{Theorem}[section]
\newtheorem{lemma}[theorem]{Lemma}

\newtheorem{corollary}[theorem]{Corollary}

\newtheorem{remark}[theorem]{Remark}

\theoremstyle{definition}

\renewcommand{\tilde}{\widetilde}          
\DeclareMathSymbol{\leqslant}{\mathalpha}{AMSa}{"36} 
\DeclareMathSymbol{\geqslant}{\mathalpha}{AMSa}{"3E} 
\DeclareMathSymbol{\eset}{\mathalpha}{AMSb}{"3F}     
\renewcommand{\leq}{\;\leqslant\;}                   
\renewcommand{\geq}{\;\geqslant\;}                   

\newcommand{\R}{\mathbb{R}}

\newcommand{\E}{\mathds{E}}
\newcommand{\Pb}{\mathds{P}}
\newcommand{\ind}{\mathds{1}}

\def\H{\mathbb{H}}

\def\bi{\begin{itemize}}
\def\ei{\end{itemize}}

\def\bnum{\begin{enumerate}}
\def\enum{\end{enumerate}}
\def\LB{\mathcal{B}} 
\def\<#1{\langle #1 \rangle}

\def\p{\mathbf{p}}

\title{Spectral dimension of  Liouville quantum gravity}
\author{R\'emi Rhodes
 \footnote{
 Universit{\'e} Paris-Dauphine, Ceremade, F-75016 Paris, France. Partially supported by grant ANR-11-JCJC  CHAMU
 }
 \and 
 Vincent Vargas  \footnote{
 Universit{\'e} Paris-Dauphine, Ceremade, F-75016 Paris, France. Partially supported by grant ANR-11-JCJC  CHAMU
 }
 }

\date{}

\begin{document}

\maketitle
 
\begin{abstract} This paper is concerned with computing the spectral dimension of $2d$-Liouville quantum gravity. As a warm-up, we first treat the simple case of boundary Liouville quantum gravity. We prove that the spectral dimension is $1$ via an exact expression for the boundary Liouville Brownian motion and heat kernel. Then we treat the $2d$-case via a decomposition of time integral transforms of the Liouville heat kernel into Gaussian multiplicative chaos of Brownian bridges. We show that the spectral dimension is $2$ in this case, as announced by physicists (see Ambj\o rn and al. in \cite{amb})  fifteen years ago.
\end{abstract}
\vspace{0.5cm}
\footnotesize


\noindent{\bf Key words or phrases:} Liouville quantum gravity, Liouville Brownian motion, Gaussian multiplicative chaos, Liouville heat kernel, spectral dimension.

\normalsize

%
%
%
\tableofcontents

\section{Introduction}
String theory is an attempt to overcome the difficulties encountered in the quantization of $4d$ gravity by replacing particles by string-like one dimensional objects, which describe some two dimensional "worldsheet" $\Sigma$ as they evolve in time. Polyakov \cite{Pol} showed that such theories could be interpreted as theories of two dimensional quantum gravity, in which the worldsheet is exchanged with the space time and the string coordinate $H(\sigma)$ ($\sigma\in\Sigma$) is considered as a $d$-dimensional matter field defined on $\Sigma$. The worldsheet of the string can be seen as a two dimensional random surface  embedded in $d$-dimensional space and equipped with a metric  $g$. This metric is a random variable, which takes on the form \cite{Pol,cf:KPZ,cf:Da} (we consider an Euclidean background metric for simplicity):
$$g(z)=e^{\gamma X(z)}dz^2,$$
where $\gamma$ is a coupling constant expressed in terms of the central charge of the matter field $H$ and $X$ is a random field, the fluctuations of which   are governed by the Liouville action. In critical Liouville quantum gravity, this action turns the field $X$ into a Free Field, with appropriate boundary conditions. The reader is referred to \cite{david-hd,cf:Da,DFGZ,DistKa,bourbaki,GM,glimm,cf:KPZ,Nak,Pol} for more insights on the subject.

Among the main objectives of $2d$-quantum gravity is the understanding of the fractal structure of space-time. In the case of pure gravity, i.e. when the central charge of the conformal matter is $c=0$ yielding $\gamma=\sqrt{8/3}$, it is known that the intrinsic Hausdorff dimension is $4$ whereas the intrinsic Hausdorff dimension related to the background metric is $2$ (see \cite{KKMW,AW} on the physics side and \cite{LeGallS,legall,miermontsphere,miermont} for more recent mathematical references). The Knizhnik-Polyakov-Zamolodchikov formula \cite{cf:KPZ} relating the scaling exponents of the matter field computed under the background metric to those computed under the metric $g$ is another striking feature of the fractal structure of space-time, which has recently received growing interest \cite{Rnew4,Benj,David-KPZ,Rnew12,cf:DuSh,Rnew10}. Another interesting feature of the fractal structure of space-time is the notion of spectral dimension. It has been proved in \cite{GRV,GRV-FD} that one can associate to the Liouville metric a Brownian motion and a heat kernel $\p_t(x,y)$. The spectral dimension is about the short time behavior of the heat kernel along the diagonal. One possible way to rigorously define the spectral dimension is:
\begin{equation}\label{spec:exact}
d_S=2\lim_{t\to 0}\frac{\ln \p_t(x,x)}{-\ln t}.
\end{equation}
In a loose sense, this means that:
\begin{equation}\label{spec:sim}
\p_t(x,x)\simeq \frac{C_x}{t^{d_S/2}},\quad \text{as }t\to 0.
\end{equation}

For instance, the spectral dimension of the Euclidean $\R^d$-space is $d_S=d$.  After averaging in some sense over geometries (the field $X$) and points $x$, the authors in \cite{amb} have heuristically obtained fifteen years ago that the spectral dimension of $2d$-quantum gravity is $2$, irrespectively of the coupling constant $\gamma$. From a rigorous angle, giving sense to the above limit is a difficult task because, for the time being, it is not clear that the Liouville heat kernel  possesses any kind of regularity, apart from measurability. Yet, we want to reinforce the statement made in \cite{amb} in a quenched version in $X$ and everywhere in $x$. The price to pay is that, instead of considering the heat kernel on the diagonal as a function, we will consider it as a measure $\p_t(x,x)\,dt$ on the positive reals. Then we focus on the mass of this measure at the neighborhood of $t=0$. Our argument relies on the fact that the mass at $0$ of this measure  can be identified via the behavior of its Laplace transform. A similar idea has already been used by M. Bauer and F. David  in  \cite{David-KPZ} where a heat kernel based KPZ relation is established via the Mellin-Barnes transform of the heat kernel. Then we will use the fact that Laplace transforms of the Liouville heat kernel are perfectly defined for all $x$  on the diagonal (and are even a continuous function of $x$) and may be decomposed into functionals of Gaussian multiplicative chaos along Brownian bridges. This allows us to have a complete description of the Laplace transforms of the Liouville heat kernel and prove  that the spectral dimension is $2$ almost surely in $X$ for every $x\in\R^2$.
 
We will also consider boundary Liouville quantum gravity, corresponding to the description of open strings. In that case, the worldsheet $\Sigma$ possesses a boundary and the field $ X$ fluctuates along the boundary.  We will treat the case when the worldsheet $\Sigma$ is the half-plane. In critical Liouville quantum gravity, the field $X$ is then a free field on the half-plane with free boundary conditions. We are thus led to considering a metric on the real line of the form
$$g(z)=e^{\gamma X(z)}dz^2,$$
 where $X$ is a centered Gaussian distribution with short scale logarithmic correlations. We first give a complete description of boundary Liouville quantum gravity by expliciting its metric, heat kernel and Brownian motion. In particular, the spectral dimension is $1$ and thus coincides with its intrinsic Hausdorff dimension. Though straightforward, it seems that such explicit formulas have never been written. 
 
Finally, we just mention that the partition function of non-critical $2d$-quantum gravity is absolutely continuous with respect to critical $2d$-quantum gravity so that its spectral dimension is also $2$.
\section{Boundary Liouville quantum gravity}\label{sec.1d}
We begin with the case of boundary quantum gravity since it provides an interesting intuition for the $2d$-case. However, as soon as the problem has been stated properly, computations are straightforward from the mathematical angle. So we just explain the main lines without giving details. This section may just be seen as a warm-up, yet it is instructive.
 
\subsection{Reminder about $1d$-Riemannian structures}\label{sec:riem}
Let $g(x)dx^2$ be a smooth metric tensor on $\R$, with $g>0$. We denote by $\triangle_g$ the associated Laplace-Beltrami operator, i.e.:
$$\triangle_g=g^{-1}(x)\triangle+g^{-1/2 }(x)\partial_x (g^{-1/2 }(\cdot))(x)\partial_x.$$
We define the smooth strictly increasing function on $\R$:
$$\varphi(x)=\int_0^xg^{1/2 }(x) \,du.$$
If we further assume that $\varphi$ maps $\R$ onto $\R$ then $\varphi$ is a Riemannian diffeomorphism between $(\R,g)$ and $(\R,1)$ where $1$ stands for the Euclidean metric. It is then obvious to describe the geometric objects associated to $g$. The Riemannian distance is  $d(x,y)=|\varphi(x)-\varphi(y)|$. The Brownian motion associated to $g$ is the solution of the SDE:
\begin{equation}
\LB^x_t=x+\int_0^t\frac{1}{2}g^{-1/2 }(\LB^x_r)\partial_x (g^{-1/2 })(\LB^x_r)\,dr+\int_0^tg^{-1/2 }(\LB^x_r)\,dB_r,
\end{equation}
where $B$ is a standard Brownian motion on $\R$. It is plain to see that: 
\begin{equation}
\forall t\geq 0,\quad \LB^x_t=\varphi^{-1}(\varphi(x)+B_t).
\end{equation}
This can be seen because diffeormorphisms of Riemannian manifolds maps Brownian motion to  Brownian motion (also follows from   $\triangle_g \varphi=0$ and $\partial_g\varphi=1$). 

In what follows, we want to replace the function $g$ by the exponential of a free field, and more generally by a centered Gaussian distribution with short scale logarithmic correlations. So, we consider a centered Gaussian distribution $X$ on $\R$ with covariance kernel of the type
$$K(x,y)=\ln \frac{1}{|x-y|}+g(x,y),$$ for some continuous and bounded function $g$ on $\R\times\R$.
For instance, one may consider a massive free field $X$ with free boundary condition on the half-plane $\H\subset \R^2$ and consider the trace of the field along the real line. 

\subsection{Standard phase}
We consider $\gamma<2\sqrt{2}$. Let us stress that, in the physics literature, people often use the coupling constant $\gamma'=\frac{\gamma}{\sqrt{2}}\in [0,2[$ instead of $\gamma$.  We want to investigate the random metric on $\R$ formally defined by $e^{\gamma X(x)}dx^2$. Following subsection \ref{sec:riem}, one wishes to define
 the random mapping: 
$$\varphi(x)=\int_0^x e^{\frac{\gamma}{2} X(r)-\frac{\gamma^2}{8}\E[X^2]}dr.$$
This can be done via cutoff approximations of the field $X$ as prescribed by standard theory of Gaussian multiplicative chaos (see \cite{cf:Kah,review} for instance). Almost surely in the field $X$, this mapping is continuous, strictly increasing and satisfies
$$\lim_{x\pm\infty} \varphi(x)=\pm \infty.$$
It can be seen as an isomorphism between $\R$ equipped with the Euclidean metric and $\R$ equipped with the metric  $g=e^{\gamma X(x)-\frac{\gamma^2}{4}\E[X^2]}dx^2$. The boundary quantum metric is given by:
$$d(x,y)=|\varphi(x)-\varphi(y)|.$$
The volume form matches:
$$\forall A \in\mathcal{B}(\R),\quad M(A)=\int_A e^{\frac{\gamma}{2} X(r)-\frac{\gamma^2}{8}\E[X^2]}dr$$
and the boundary Liouville Brownian  motion $\LB^x$, starting from $x$, reads:
$$\LB^x_t=\varphi^{-1}(\varphi(x)+B_t).$$
It is defined almost surely in $X$ for all starting points $x\in\R$. It is further a strong Feller Markov process with continuous sample paths. Observe that this process is formally the solution of the SDE:
\begin{equation}
\LB^x_t=x-\frac{\gamma^2}{4}\int_0^t\partial_xX(\LB^x_r) e^{- \gamma X(\LB^x_r)+\frac{\gamma^2}{4}\E[X^2] } \,dr+\int_0^te^{- \frac{\gamma}{2}  X(\LB^x_r)+\frac{\gamma^2}{8}\E[X^2] } \,dB_r.
\end{equation}
The boundary Brownian motion admits a heat kernel. More precisely, there exists  a continuous family $\p_t(\cdot,\cdot))_t$ such that for all $x,y,t$:
$$\E[f(\LB^x_t)]=\int_\R \p_t(x,y)f(y)\,M(dy).$$ By making straightforward changes of variables, one can obtain the following explicit form for the heat kernel:
$$\p_t(x,y)=\frac{1}{\sqrt{2\pi t}}e^{-\frac{d(x,y)^2}{2t}}.$$

\begin{corollary}
The spectral dimension of boundary Liouville quantum gravity, defined by:
$$d_S=-2\lim_{t\to 0}\frac{\ln \p_t(x,x)}{\ln t} $$  is given by $d_S=1$.
\end{corollary}
Observe that it is here obvious to see that the intrinsic Hausdorff dimension of boundary Liouville quantum gravity is also $1$. Let us further stress that the topology of the quantum metric is Euclidean: this results from the fact that the field $\varphi$ is continuous and strictly increasing. Yet the quantum metric is not equivalent to the Euclidean one: this results from the multifractal analysis of the field $\varphi$, i.e. a study of the decrease of the size of balls (see \cite{BS,review}).  
\subsection{Critical phase}\label{sec.crit}

The critical phase corresponds to  the case  $\gamma=2\sqrt{2}$.  The mathematical proofs do not differ from the subcritical situation so that we only summarize the results. 
By considering a family of smooth cut off approximations $(X_\epsilon)_\epsilon$ of the field $X$, we introduce the random smooth metric tensor:
$$g_\epsilon=e^{2\sqrt{2} X_\epsilon(x)-2\E[X_\epsilon^2]}dx^2$$
and we define the random mapping: 
$$\varphi_\epsilon(x)=\int_0^x e^{\sqrt{2}X_\epsilon(r)-\E[X_\epsilon^2]}dr.$$
It is proved in \cite{Rnew7,Rnew12} that, almost surely in $X$, the following convergence: 
$$\sqrt{-\ln\epsilon}\varphi_\epsilon\to \varphi(x)=\sqrt{2/\pi}\int_0^x (\sqrt{2}\E[X^2]-X(r))e^{\sqrt{2} X(r)-\E[X^2]}dr$$ holds in the space $C(\R)$. The limiting mapping will be called the critical diffeomorphism. The mapping $\varphi$ is continuous ans strictly increasing and goes to $\pm \infty$ as $x\to\pm \infty$. Once again, the metric matches:
$$d(x,y)=|\varphi(x)-\varphi(y)|$$ 
and the volume form  $$\forall A \in\mathcal{B}(\R),\quad M'(A)= \sqrt{2/\pi}\int_A (\sqrt{2}\E[X^2]-X(r))e^{\sqrt{2} X(r)-\E[X^2]}dr.$$

Almost surely in $X$ and $B$, for all $x\in\R$, the critical boundary Liouville Brownian motion is given by:
$$\LB^x_t=\varphi^{-1}(\varphi(x)+B_t).$$
It is a Feller Markov process with continuous sample paths. It admits a heat kernel of the same form as in the  standard phase. Therefore the spectral dimension of critical boundary Liouville quantum gravity is still $1$, as well as the intrinsic  Hausdorff dimension.

\section{$2d$-Liouville quantum gravity}
In this section, we consider the same setup as in \cite{GRV,GRV-FD}. We just outline the main tools needed in this paper and the reader is referred to \cite{GRV,GRV-FD} for further details. We consider a whole plane Massive Gaussian Free Field (MFF) $X$ with mass $m$ (see \cite{glimm,She07} for an overview of the construction of the MFF and applications). Its covariance kernel is thus given by 
\begin{equation}\label{MFF1}
\forall x,y \in \R^2,\quad G_m(x,y)=\int_0^{\infty}e^{-\frac{m^2}{2}u-\frac{|x-y|^2}{2u}}\frac{du}{2 u}=\ln_+\frac{1}{|x-y|}+g_m(x,y).
\end{equation}
for some continuous and bounded function $g_m$, which decays exponentially fast to $0$ when $|x-y|\to\infty$ (recall that $\ln_+(x)= \max(0,\ln x)$ for $x>0$). 

We consider a coupling constant $\gamma \in [0,2[$ and consider the formal metric tensor:
$$g=e^{\gamma X(x)-\frac{\gamma^2}{2}\E[X^2]}\,dx^2.$$
This metric tensor has a volume form, which is nothing but a Gaussian multiplicative chaos \cite{cf:Kah,review} with respect to the Lebesgue measure $dx$:
\begin{equation}
M(dx)=e^{ \gamma X(x)-\frac{\gamma^2}{2}\E[X(x)^2]}\,dx.
\end{equation}
One may also associate to this metric a Brownian motion, called Liouville Brownian motion (LBM for short). More precisely, almost surely in $X$, for all $x\in\R^2$, the law of the LBM starting from $x$  is given by
$$ \LB^x_t=B^x_{F(x,t)^{-1}}$$ where $B^x$ is a standard two dimensional Brownian motion starting from $x$ and the random mapping $F$ is defined by
$$F(x,t)=\int_0^t e^{ \gamma X(B^x_r)-\frac{\gamma^2}{2}\E[X^2]}\,dr.$$
It is a Feller Markov process with continuous sample paths. It is also shown in \cite{GRV-FD} that, almost surely in $X$ and for all starting point $x\in\R^2$, the law of the LBM is absolutely continuous with respect to the Liouville measure $M$, thus giving the existence of the Liouville heat kernel $\p_t(x,y)$ such that for all $x\in\R^2$ and all measurable bounded function $f$:
\begin{equation}
\E[f(\LB^x_t)]=\int_{\R^2}f(y)\p_t(x,y)\,M(dy).
\end{equation}
In what follows, we will also consider the heat kernel $p_t(x,y)$  of the standard Brownian motion on $\R^2$.

\subsection{Brownian bridge decomposition}
Let us denote by $(B^{x,y,t}_s)_{0\leq s \leq t}$ a (standard or Euclidean) Brownian bridge between $x$ and $y$ with lifetime $t$. We start with a standard lemma for Brownian bridges in dimension 2:
\begin{lemma}\label{bridge}
Let $(B_u^x)_{u \leq t}$ be a Brownian motion starting from $x$. We have the following absolute continuity relation for the Brownian bridge when $s<t$:
\begin{equation}\label{eq:bridge}
\E[  F( (B^{x,y,t}_u)_{u \leq s}   ) ]  = \E^x[   F( (B_u^x)_{u \leq s}   )  \frac{t}{t-s}e^{\frac{|y-x|^2}{2 t}-\frac{|B_s^x-y|^2}{2 (t-s)}  }     ]
\end{equation}
\end{lemma}

Here and in the sequel, $\E$ means that we take expectation with respect to the law of the Brownian bridge $B^{x,y,t}$, whereas $\E^x$ means expectation with respect to a Brownian motion started at $x$. We will also denote by $\E$ (or $\E^X$ in case there could be an ambiguity) expectation with respect to the free field $X$.

For each $x\in\R^2$ and $s\in [0,t]$, we define
\begin{equation}
 F(x,y,t,s)=\int_0^se^{\gamma X(B^{x,y,t}_r)-\frac{\gamma^2}{2}\E[X^2]}\,dr.  
\end{equation}
Actually, a rigorous definition of such an object is not straightforward. One has to introduce a cut off approximation sequence $(X_n)_n$ of the field $X$ as considered in \cite{GRV} and define $F(x,y,t,\cdot)$ as the limit in law in the space $C([0,t])$ of the sequence:
$$F_n(x,y,t,s)=\int_0^se^{\gamma X_n(B^{x,y,t}_r)-\frac{\gamma^2}{2}\E[X_n^2]}\,dr.$$  
We state the two following theorems:

\begin{theorem}
Assume $\gamma<2$. For each $x,y\in\R^2$ and $t\geq 0$, almost surely in $X$ and in $B^{x,y,t}$, the family of random mappings  $(F_n(x,y,t,\cdot))_n$ converges in $C([0,t],\R_+)$ towards a limiting strictly increasing continuous mapping $F(x,y,t,\cdot)$.
\end{theorem}

\begin{theorem}\label{thtouspoints}
Assume $\gamma<2$. Almost surely in $X$, for all $x,y\in\R^2$ and $t\geq 0$, the law under $\Pb^{B^{x,y,t}}$ of the random mappings $(F_n(x,y,t,\cdot))_n$ converges in $C([0,t],\R_+)$ towards the law of a random mapping,  still denoted by $F(x,y,t,\cdot)$. Under $\Pb^{B^{x,y,t}}$, the mapping $F(x,y,t,\cdot)$ is strictly increasing and continuous. Furthermore, for each bounded continuous function $K:C([0,t],\R_+)\to\R$, the mapping 
$$(x,y)\mapsto \E^{B^{x,y,t}}[K(F(x,y,t,\cdot))]$$ is continuous. Finally, 
$$\E^x[   G( F(x,t) )   |  B_t=y  ]=\E\big[G(F(x,y,t,t))\big]$$ for all nonnegative measurable function $G$. 
\end{theorem}

To prove the above two theorems, a rigorous argument just boils down to reproducing the arguments of \cite{GRV}. In particular, to prove theorem \ref{thtouspoints},  one can adapt the coupling technique between two Brownian motions starting from distinct points introduced in \cite{GRV} (or see lemma \ref{lem:coupling} below). The fact that this technique can be used here in the context of bridges is a consequence of the absolute continuity formula (\ref{eq:bridge}) and the time reversal symmetry of Brownian bridges.


We then prove the following Brownian bridge decomposition of integral transforms of the Liouville heat kernel:
\begin{theorem}\label{generalformula}
Let $G:\R_+\to\R_+$ be a   continuous function. Then for all $x\in\R^2$ and $M$-almost every $y\in\R^2$:
\begin{equation} \label{transformBB}
\int_{0}^\infty  G(t)    \p_t(x,y) dt =   \int_0^\infty  \E[ G(  F(x,y,t,t))  ]  \frac{e^{-\frac{|y-x|^2}{2t}}}{2 \pi t}  dt  
\end{equation}
\end{theorem}
 
\noindent {\it Proof.}  We have for all bounded continuous function $f$ on $\R^2$:
\begin{align*}
 \int_{\R^2}\Big( \int_{0}^\infty  G(t) &  \p_t(x,y) dt \Big)f(y) M(dy)\\ = &\int_{0}^\infty  G(t)\Big(\int_{\R^2}   \p_t(x,y)f(y)  M(dy)\Big)\,dt\\
 =& \E^x\Big[\int_{0}^\infty  G(t)f(B^x_{F(x,t)^{-1}})\,dt\Big]\\
  =& \E^x\Big[\int_{0}^\infty  G(F(x,t))f(B^x_{t})\,F(x,dt)\Big]\\
    =& \E^x\Big[\int_{0}^\infty  \E\big[G\big(F(x,y,t,t)\big)\big]_{y=B^x_{t}}f(B^x_{t})\,F(x,dt)\Big]\\
     =& \int_{\R^2} \int_{0}^\infty  \E\big[G\big(F(x,y,t,t)\big)\big]f(y)p_t(x,y)\,dt\,M(dy) 
\end{align*} 
We deduce that, almost surely in $X$, for all $x\in\R^2$ and for $M$-almost every $y$:
$$\int_{0}^\infty  G(t)     \p_t(x,y) dt= \int_{0}^\infty  G\big(F(x,y,t,t)\big)p_t(x,y)\,dt.$$
This completes the proof.\qed

\begin{remark}
Note that theorem \ref{generalformula} is quite general and in particular gives a useful and rather explicit formula for the resolvent density (for $G(t)=e^{-\lambda t}$ with $\lambda>0$) or the Mellin-Barnes transform considered in \cite{David-KPZ} (this corresponds to $G(t)=t^{s-1}$ with $s \in ]0,1[$). Furthermore, the forthcoming proofs also show that these transforms are continuous functions of $x,y$ and of the parameter $\lambda$ or $s$, depending on the considered integral transform. A more precise study of the Mellin-Barnes transform will be presented in a forthcoming work. 
\end{remark}

Now we take $G(t)=t^{\alpha}e^{-\lambda t}$ for $t\geq 0$ and $\alpha\geq 0$ to get:
\begin{equation} \label{transformBB}
\int_{0}^\infty  e^{-\lambda t}  t^{\alpha}  \p_t(x,y) dt =   \int_0^\infty  \E[  F(x,y,t,t)^\alpha e^{-\lambda F(x,y,t,t)  }]  \frac{e^{-\frac{|y-x|^2}{2t}}}{2 \pi t}  dt  .
\end{equation}
This identity is valid for every $x\in\R^2$ and $M$-almost every $y\in\R^2$. In the following, we will be interested in the diagonal behavior ($x=y$) of this quantity. We will see that one can make sense of these quantities on the diagonal.

\subsection{Spectral dimension}

We can now state our main result on the spectral dimension:

\begin{theorem}
For $\alpha\geq 0$, the right-hand side of \eqref{transformBB} admits a limit as $y\to x$ for all $x\in\R^2$.  Therefore \eqref{transformBB} makes sense on the diagonal $y=x$. Furthermore, for $\alpha>0$ and for all $x\in\R^2$:
$$\forall\lambda>0,\quad\int_{0}^\infty  e^{-\lambda t}  t^{\alpha}  \p_t(x,x) dt <+\infty.$$
For $\alpha=0$, we have for all $x\in\R^2$:
$$\forall\lambda>0,\quad \int_{0}^\infty  e^{-\lambda t}    \p_t(x,x) dt =+\infty.$$
Therefore, the spectral dimension of Liouville quantum gravity is $2$.
\end{theorem}

\begin{remark}
A stronger statement consists in proving that $$\lim_{t\to 0}t\p_t(x,x)=c_x$$ for some random constant $c_x$ such that $c_x\ne 0$ for all $x\in\R^2$. In a way, the above theorem implies that if one can show that
$$\lim_{t\to 0}t^{\nu_x/2}\p_t(x,x)=c_x$$
for some random constant $\nu_x$, then $\nu_x$ can be nothing but $2$.
\end{remark}
\vspace{2mm}

\noindent {\it Proof.} Let us begin with the case $\alpha=0$. For $\alpha=0$, we will understand the "on-diagonal" relation \eqref{transformBB} as a limit as $|x-y|\to 0$, thus being left with proving that such a limit exists. For $x,y\in\R^2$, we have:
\begin{align*}
\int_0^\infty  \E[   e^{-\lambda F(x,y,t,t)  }]  &\frac{e^{-\frac{|y-x|^2}{2t}}}{2 \pi t}  dt \\
\geq & \int_0^\infty  \E[   e^{-\lambda F(x,y,t,t)  }\ind_{\{F(x,y,t,t)\leq 1\}}]  \frac{e^{-\frac{|y-x|^2}{2t}}}{2 \pi t}  dt\\
\geq &  e^{-\lambda    }\int_0^\infty  \E[  \ind_{\{F(x,y,t,t)\leq 1\}}]  \frac{e^{-\frac{|y-x|^2}{2t}}}{2 \pi t}  dt\\
\end{align*}
Therefore, by Fatou's lemma and for $\delta$ arbitrarily small:   
\begin{align*}
\liminf_{|x-y|\to 0}\int_0^\infty  \E[   e^{-\lambda F(x,y,t,t)  }]  &\frac{e^{-\frac{|y-x|^2}{2t}}}{2 \pi t}  dt \geq     e^{-\lambda    }\int_0^\infty  \E[  \ind_{\{F(x,x,t,t)\leq 1-\delta\}}]  \frac{1}{2 \pi t}  dt.
\end{align*}
Here we have used the fact that the mapping $(x,y)\mapsto F(x,y,t,t)$ is continuous in law. Since the function $x\mapsto \ind_{\{x\leq 1\}}$ is not continuous, we may estimate from below this function with a continuous piecewise linear function $\varphi$ such that $\varphi(x)=1$ for $x\leq 1-\delta$ and $\varphi(x)=0$ for $x\geq 1$, thus the presence of a small $\delta$ in the above relation. By noticing that, for all $x\in\R^2$, $\lim_{t\to 0}\E[  \ind_{\{F(x,x,t,t)\leq 1-\delta\}}] =1$, we deduce that for all $x$: 
$$\int_0^\infty  \E[  \ind_{\{F(x,x,t,t)\leq 1-\delta\}}]  \frac{1}{2 \pi t}  dt=+\infty.$$
The proof of the case $\alpha=0$ follows.

Now we have to treat the case $\alpha>0$. If we can prove that $\E[F(x,y,t,t)^\alpha]\leq Ct^{\epsilon}$ for some $\epsilon>0$ and all $y$ in a neighborhood of $x$ and $t\leq 1$, then the part in the right-hand side of \eqref{transformBB} corresponding to: 
$$\int_0^1  \E[  F(x,y,t,t)^\alpha e^{-\lambda F(x,y,t,t)  }]  \frac{e^{-\frac{|y-x|^2}{2t}}}{2 \pi t}  dt$$
is continuous on the diagonal and:
\begin{equation}\label{est1}
\int_0^1  \E[  F(x,x,t,t)^\alpha e^{-\lambda F(x,x,t,t)  }]  \frac{1}{2 \pi t}  dt<+\infty.
\end{equation}
So let us investigate the quantity $\E[F(x,y,t,t)^\alpha]$. By using the time reversal symmetry of the Brownian bridge and the sub-additivity of the mapping $x\mapsto x^\alpha$, we get:
$$\E[F(x,y,t,t)^\alpha]\leq \E[F(x,y,t,t/2)^\alpha]+\E[F(y,x,t,t/2)^\alpha].$$
Therefore, it suffices to investigate $\E[F(y,x,t,t/2)^\alpha]$ (or equivalently $\E[F(x,y,t,t/2)^\alpha]$). We will use the fact that the law of the Brownian bridge on $[0,t/2]$ looks like the Brownian motion (see lemma \ref{bridge}). 

Indeed, using lemma \ref{bridge}, we deduce that, for some constant $C$ which does not depend on $t$ or $|y-x|\leq 1$, we have:
\begin{align*}
\E[F(y,x,t,t/2)^\alpha]e^{-\frac{|y-x|^2}{2t}}  &\leq C \E^y[F(y,t/2)]^\alpha\\
&\leq C\Big(\E^y\Big[\int_0^{t/2}e^{\gamma X(B_s^y)-\frac{\gamma^2}{2}\E[X^2]}\,ds\Big] \Big)^\alpha\\
&=C\Big(\int_0^{t/2}\int_{\R^2} p_{s}\big(y,z\big)\,M(dz)\,ds \Big)^\alpha.
\end{align*}
We state the following two lemmas which we will prove after:
\begin{lemma}\label{decomheat}
We have:
$$ \int_0^{t/2} p_{s}\big(y,z\big)\,ds\leq C(1+\ln\frac{t^{1/2}}{|y-z|})\ind_{\{|y-z|\leq t^{1/2}\}}+Ce^{-\frac{|z-y|^2}{t}}\ind_{\{|y-z|\geq t^{1/2}\}}.$$
\end{lemma}
and also:
\begin{lemma} \label{app:multform}
Let $\epsilon>0$ and $R>0$. We set $\beta=2(1-\frac{\gamma}{2})^2>0$. Almost surely in $X$, there exists a random constant $C>0$ such that:
\begin{equation*}
\sup_{x \in [-R,R]^2} M( B(x,r)  )  \leq C r^{\beta-\epsilon}, \quad \forall r\in]0,1[.
\end{equation*}
\end{lemma} 

\vspace{2mm} Now choose $R>|x|+2$ and $\epsilon=\beta/2$. We use these two lemmas to get for $|y-x|\leq 1$ and $t\leq 1$: 
 \begin{align*}
 \int_{\R^2}(1+\ln\frac{t^{1/2}}{|y-z|})&\ind_{\{|y-z|\leq t^{1/2}\}}M(dz)\\
 \leq &\sum_{n\geq 0} \int_{\R^2}(1+\ln\frac{t^{1/2}}{|y-z|}) \ind_{\{ t^{1/2}2^{-n}\leq |y-z|\leq t^{1/2}2^{-n+1}\}}M(dz)\\
 \leq & \sum_{n\geq 0} (1+n\ln 2) M(B(y, t^{1/2}2^{-n+1})\\
 \leq & t^{\frac{\beta}{4}}C\sum_{n\geq 0} (1+n\ln 2) 2^{-n\beta/2}.
\end{align*}

Now we focus on the area $|y-z|\geq t^{1/2}$. Let $\delta\in]0,1/2[$. We have:
 \begin{align*}
 \int_{\R^2}e^{-\frac{|z-y|^2}{t}}&\ind_{\{|y-z|\geq t^{1/2}\}}M(dz)\\
 \leq & M(B(y,t^{1/2-\delta}))+ \int_{\{|y-z|\geq t^{1/2-\delta},|z-x|\leq 4\}}e^{-\frac{|z-y|^2}{t}} M(dz)     \\
 &+      \int_{\{|z-x|\geq 4\}}e^{-\frac{|z-y|^2}{t}} M(dz)\\
 \leq & Ct^{(1/2-\delta)\beta/2}+e^{-t^{-2\delta}}M(B(x,4))+  \int_{\{|z-x|\geq 4\}}e^{-\frac{|z-x|^2}{2t}}\,M(dz)\\
 \leq & Ct^{(1/2-\delta)\beta/2}+e^{-t^{-2\delta}}M(B(x,4))+ C t^2 \int_{\R^2}e^{-\frac{|z-x|^2}{2}}\,M(dz)
\end{align*}
where the constant $C$ does not depend on $t\leq 1$. Almost surely in $X$, the right-hand side of the above inequality is finite and is less than $C t^\xi$ for some $\xi>0$, for all $t\leq 1$ and for some random constant $C$. By gathering all the above considerations, we deduce:
$$\E[F(x,y,t,t)^\alpha]\leq Ct^{\alpha \xi}$$ for some random constant $C$ which is finite for $|y-x|\leq 1$ and $t\leq 1$.

It remains to prove that the quantity:
\begin{equation}\label{est2}
\int_{1}^\infty   \E[  F(x,y,t,t)^\alpha e^{-\lambda F(x,y,t,t)  }]  \frac{e^{-\frac{|y-x|^2}{2t}}}{2 \pi t}  dt 
\end{equation}
 is continuous  and finite on the diagonal. Once again, we can use the absolute continuity of the law of the Brownian bridge for $|y-x|\leq 1$ to see that it suffices to prove: 
$$\int_{1}^\infty   \sup_{|y-x|\leq 1}   \E^y[  F(y,t/2)^\alpha e^{-\lambda F(y,t/2)  }]  \frac{1}{2 \pi t}  dt<+\infty.$$ By noticing that $u^\alpha e^{- \lambda u }\leq C$ for $u\geq 0$ and some constant $C$, we deduce that it suffices to prove: 
\begin{equation}\label{suptam}
\int_{1}^\infty  \sup_{|y-x|\leq 1}  \E^y[    e^{-\lambda F(y,t/2)  }]  \frac{1}{2 \pi t}  dt<+\infty.
\end{equation}
 
  \hspace{0.5 cm}
 
 \noindent
 \emph{Step one:} we first prove that $\int_{1}^\infty   \E^x[    e^{-\lambda F(x,t/2)  }]  \frac{1}{2 \pi t}  dt<+\infty$.

  \hspace{0.3 cm}

 Without loss of generality, we may assume $x=0$. In this first step, for the sake of clarity, we skip the dependency on the starting point $x$ in $F(x,t)$ and simply write $F(t)$ for $F(0,t)$. We will also write $F(s,t)$ for $F(t)-F(s)$.
 
We first assume that $X$ is decorrelated at distance $1$, meaning that $\E[X(x)X(y)]=0$ for $|x-y|\geq 1$ (of course, this formal relation $\E[X(x)X(y)]=0$ is not rigorous as the field $X$ does not makes sense pointwise  but it is straightforward to see how to make sense of it). We introduce the following family of increasing stopping times $(T_n, \tilde{T}_n)$ for $n\geq 1$:
 \begin{equation*}
 T_n= \inf \lbrace t \geq \tilde{T}_{n-1}, \: |B_t|= 3 n    \rbrace, \; \; \tilde{T}_n= \inf \lbrace t \geq T_{n}, \: |B_t-B_{T_n} | \geq 1    \rbrace, \; \; 
 \end{equation*}
 where by convention $T_0=\tilde{T}_0=0$. We set $N(t)= \sup \lbrace n,; \;  \tilde{T}_n \leq t \rbrace$. Now we have that:
 \begin{align*} 
  \E^X \E^0[    e^{-\lambda F(0,t)  }]  & \leq   \E^X \E^0[   \Pi_{n=1}^{N(t)}   e^{-\lambda F( T_n,  \tilde{T}_n)} ] \\
& =  \E^0[   \Pi_{n=1}^{N(t)} \E^X[  e^{-\lambda F( T_n,  \tilde{T}_n)} ] ]  \\
& =  \E^0[   \Pi_{n=1}^{N(t)} \E^X[  e^{-\lambda \int_{T_n}^{\tilde{T}_n} e^{\gamma X(B_t-B_{ T_n})-\frac{\gamma^2}{2}\E[ X^2  ] }\,dt} ] ] .
 \end{align*}
We will exploit the fact that the Brownian motions $(B_t-B_{ T_n})_{ t \in [T_n, \tilde{T}_n]}$ are independent (note however that they are not independent of $N(t)$). Now fix $\epsilon >0$. We have:
 \begin{align*} 
& =  \E^0[   \Pi_{n=1}^{N(t)} \E^X[  e^{-\lambda \int_{T_n}^{\tilde{T}_n} e^{\gamma X(B_t-B_{ T_n})-\frac{\gamma^2}{2}\E[ X^2  ] }\,dt} ] ]  \\
& \leq \Pb^0(N(t) \leq t^{\frac{1}{2}-\epsilon}) +  \E[ 1_{N(t) \geq t^{\frac{1}{2}-\epsilon}}  \Pi_{n=1}^{t^{\frac{1}{2}-\epsilon}} \E^X[  e^{-\lambda \int_{T_n}^{\tilde{T}_n} e^{\gamma X(B_t-B_{ T_n})-\frac{\gamma^2}{2}\E[ X^2  ] } dt } ] ]  \\
& \leq      \Pb^0 (  \sup_{0 \leq s \leq t} |B_s| \leq t^{\frac{1}{2}-\epsilon} )  +   \E^0[  \Pi_{n=1}^{t^{\frac{1}{2}-\epsilon}} \E^X[  e^{-\lambda \int_{T_n}^{\tilde{T}_n} e^{\gamma X(B_t-B_{ T_n})-\frac{\gamma^2}{2}\E[ X^2  ] } dt } ]  ]    \\ 
&  \leq  e^{-c t^\epsilon} +    \Big(  \E^X[  \E^0[ e^{-  \lambda \int_{0}^{T} e^{\gamma X(B_t)-\frac{\gamma^2}{2}\E[ X^2  ] } dt} ]] \Big) ^{ t^{\frac{1}{2}-\epsilon} }  
\end{align*}
where $T= \inf \lbrace t \geq 0, \: |B_t|= 1 \rbrace$. We conclude with the help of the following strict inequality: $\E^X[  \E^0[ e^{-  \lambda \int_{0}^{T} e^{\gamma X(B_t)-\frac{\gamma^2}{2}\E[ X^2  ] }dt} ]]  < 1$.

We can get rid of the restriction that the field $X$ be decorrelated at distance $1$ by using  Kahane's convexity inequalities \cite{cf:Kah}. Indeed, observe that the covariance kernel of the field $X$ is then dominated by that of $X'+Y$ where $X'$ is a centered Gaussian distribution with covariance kernel given by:
$$\E[X'(x)X'(y)]=\ln_+\frac{1}{|y-x|}$$ and $Y$ is a centered Gaussian random variable with variance $\sup_{\R^2\times\R^2}|g(x,y)|$ and independent of $X'$ (and of $B$ too). The field $X'$ is thus decorrelated at distance $1$. We can then proceed 
 along the same line as above and to obtain in the end an estimate of the type:
 \begin{align*} 
  \E^X \E^0[    e^{-\lambda F(0,t)  }]  & \leq  \E^Y \left[   ( \E^{X'}\E^0[ e^{-  \lambda e^Y Z } ] )^{ t^{\frac{1}{2}-\epsilon} } \right ] 
 \end{align*}
where $Z = \int_{0}^{T} e^{\gamma X'(B_t)-\frac{\gamma^2}{2}\E[ (X')^2  ] } dt $. We introduce an i.i.d. sequence of random variables $(Z_i)_{1 \leq i \leq t^{\frac{1}{2}-\epsilon} /2 }$ with law $Z$ and we denote by $\E$ expectation with respect to this sequence. We have:
\begin{align*} 
   \E^Y \left[   ( \E^{X'}\E^0[ e^{-  \lambda e^Y Z } ] )^{ t^{\frac{1}{2}-\epsilon} } \right ]  
 & = \E [  \E^Y [    e^{-  \lambda e^Y \sum_{n=1}^{  t^{\frac{1}{2}-\epsilon} /2    }  Z_n  }  ]     ]  \\
 & \leq C \E [   \frac{1}{    \sum_{n=1}^{  t^{\frac{1}{2}-\epsilon} /2    }  Z_n  }   ]  \\
& \leq \frac{C}{t^\alpha}
 \end{align*}
 for some $\alpha >0$. This concludes step one.
 
 \hspace{0.5 cm}
 
  \noindent
 \emph{Step two:} we show that $\int_{1}^\infty  \sup_{|y-x|\leq 1}  \E^y[    e^{-\lambda F(y,t/2)  }]  \frac{1}{2 \pi t}  dt<+\infty$.
 
  \hspace{0.3 cm}

 Now, we use the following coupling lemma (this is a slight variant of the coupling lemma used in \cite{GRV}):
 \begin{lemma}\label{lem:coupling}
Fix $y_0\in\R^2$ and let us start a Brownian motion $B^{y_0}$ from $y_0$. Let us consider another independent Brownian motion $B$ starting from $0$ and denote by $B^y$, for some $y\in\R^2$, the Brownian motion $B^y=y+B$. Let us denote by $\tau^y_i$ ($i=1$ or $2$) the first time at which the $i$-th components of $B^{y_0}$ and $B^y$ coincide:
$$\tau_1^y=\inf\{u>0;B^{1,y_0}_u=B^{1,y}_u\}, \; \; \tau_2^y=\inf\{u>0;B^{2,y_0}_u=B^{2,y}_u\}$$ 
We set $\tau^{y_0,y} = \sup( \tau_1^y, \tau_2^y  )$. The random process $\overline{B}^{y_0,y}$ defined by 
$$\overline{B}^{y_0,y}_t=\left\{ 
\begin{array}{lllllll}
(B^{1,y_0}_t,B^{2,y_0}_t)  & \text{if}  &  t\leq \min( \tau_1^y, \tau_2^y)   \\
(B^{1,y}_{t}, B^{2,y_0}_t)  & \text{if}  &  \tau_1^y < t \leq \tau^{y,y_0} &  \text{or} &  (B^{1,y_0}_{t}, B^{2,y}_t)  &  \text{if} &   \tau_2^y < t \leq \tau^{y_0,y}  \\
(B^{1,y}_{t}, B^{2,y}_t)  & \text{if}  &  \tau^{y_0,y}<t .
\end{array}
\right.$$ is  a new Brownian motion on $\R^2$ starting from $y_0$, and coincides with $B^y$ for all times $t>\tau^{y_0,y}$. Furthermore, as $y\to y_0$, we have for all $\eta>0$:
$$\forall \eta>0,\quad \lim_{y\to y_0}\Pb(\tau^{y_0,y}>\eta)= 0, \quad \sup_{|y-y_0|\leq 1}\E^y[\ln \tau^{y_0,y}]<+\infty.$$
\end{lemma}

We set $\bar{x}=x+(1,1)$. Let $y \in x+[0,1]^2$. In fact, by a straightforward generalization of the above procedure, one can couple two Brownian motions $B^y$, $B^{\bar{x}}$ starting from $y$ and $\bar{x}$ to a Brownian motion $B^x$ (starting from $x$) such that $\tau^{x,y} \leq \tau^{x,\bar{x}}$. Indeed, in the above lemma, take the same driving Brownian motion for $B^y$ and $B^{\bar{x}}$. Hence, we get:
\begin{equation*}
\sup_{ y  \in x+[0,1]^2    } \E^y[    e^{-\lambda F(y,t/2)  }] \leq \E^{\bar{x}}  [ \ind_{\{\tau^{x,\bar{x}} >t/2\}} + \ind_{\{\tau^{x,\bar{x}} \leq t/2\}} \E^{B_{\tau^{x,\bar{x}}}^{\bar{x}}} [e^{-\lambda F( B_{\tau^{x,\bar{x}}}^{\bar{x}}  ,t/2-\tau^{x,\bar{x}} )  }     ]   ]. 
\end{equation*} 
 From this we deduce the following bound:
 \begin{align*}
&  \int_{1}^\infty  \sup_{|y-x|\leq 1}  \E^y[    e^{-\lambda F(y,t/2)  }]  \frac{dt}{t}     \\
 & \leq  \int_{1}^\infty  \E^{\bar{x}}  [ \ind_{\{\tau^{x,\bar{x}} >t/2\}} + \ind_{\{\tau^{x,\bar{x}} \leq t/2\}} \E^{B_{\tau^{x,\bar{x}}}^{\bar{x}}} [e^{-\lambda F( B_{\tau^{x,\bar{x}}}^{\bar{x}}  ,t/2-\tau^{x,\bar{x}} )  }     ]   ]  \frac{dt} {t}      \\
& \leq \E^{\bar{x}}  [ \ln \sup( \tau^{x,\bar{x}} ,1) ] +     \E^{\bar{x}}  [   \int_{2 \tau^{x,\bar{x}} }^\infty   \E^{B_{\tau^{x,\bar{x}}}^{\bar{x}}} [e^{-\lambda F( B_{\tau^{x,\bar{x}}}^{\bar{x}}  ,t/2-\tau^{x,\bar{x}} )  }     ]       \frac{dt}{t}  ].    
 \end{align*}
 By stationarity of the field $X$ and applying step one, we get the existence of $C,\alpha>0$ (independent from everything) such that:
 \begin{equation*}
  \E^X   \E^{B_{\tau^{x,\bar{x}}}^{\bar{x}}} [e^{-\lambda F( B_{\tau^{x,\bar{x}}}^{\bar{x}}  ,t/2-\tau^{x,\bar{x}} )  }     ]   ]  \leq \frac{C}{(t/2- \tau^{x,\bar{x}} +1)^\alpha}.
 \end{equation*}
Therefore, we get:
 \begin{align*}
 & \E^X   [  \E^{\bar{x}}  [   \int_{2 \tau^{x,\bar{x}} }^\infty   \E^{B_{\tau^{x,\bar{x}}}^{\bar{x}}} [e^{-\lambda F( B_{\tau^{x,\bar{x}}}^{\bar{x}}  ,t/2-\tau^{x,\bar{x}} )  }     ]   ]    \frac{dt}{t}  ] \\
 & =   \E^{\bar{x}}  [   \int_{2 \tau^{x,\bar{x}} }^\infty  \E^X  [  \E^{B_{\tau^{x,\bar{x}}}^{\bar{x}}} [e^{-\lambda F( B_{\tau^{x,\bar{x}}}^{\bar{x}}  ,t/2-\tau^{x,\bar{x}} )  }     ]]     \frac{dt}{t} ]  \\
 &  \leq  C \E^{\bar{x}}  [   \int_{2 \tau^{x,\bar{x}} }^\infty        \frac{dt}{t (t/2- \tau^{x,\bar{x}} +1)^\alpha}   ]   \\
 & \leq C, \\  
  \end{align*} 
 hence the desired result. 
 
\qed

  We now give the proofs of the intermediate lemmas:

\vspace{2mm}

\noindent {\it Proof of Lemma \ref{decomheat}.}   Let us set:
\begin{align*}
f_{y,t}(z)=& \int_0^{t/2} p_{s} (y,z )\,ds=\int_0^{t/2} e^{-\frac{|z-y|^2}{2u}}\,\frac{du}{2\pi u}=\int_0^{\frac{t}{2|y-z|^2}} e^{-\frac{1}{2u}}\,\frac{du}{2\pi u}.
\end{align*}
 If $|z-y|^2\leq t$, we have:
 \begin{align*}
f_{y,t}(z)\leq & \int_0^{\frac{1}{2}} e^{-\frac{1}{2u}}\,\frac{du}{2\pi u}+ \int_{\frac{1}{2}}^{\frac{t}{2|y-z|^2}} \,\frac{du}{2\pi u}\\
\leq &C\Big(1+\ln\frac{t^{1/2}}{|y-z|}\Big).
\end{align*}
If  $|z-y|^2\geq t$, we have:
 \begin{align*}
f_{y,t}(z)= &\frac{1}{4\pi } \int_{2\frac{|z-y|^2}{t}}^{+\infty} e^{-\frac{u}{2}}\,du\\
\leq &Ce^{-\frac{|z-y|^2}{t}},
\end{align*}
which completes the proof.\qed

\vspace{2mm}

\noindent {\it Proof of Lemma \ref{app:multform}.} See \cite{GRV} or also \cite{finnish} for even finer estimates. \qed


 \end{document}